\UseRawInputEncoding
\pdfoutput=1
\documentclass[a4paper,12pt]{amsart}

\usepackage{amsmath}
\usepackage{amssymb}
\usepackage{graphicx}

\theoremstyle{plain}

\theoremstyle{definition}

\newtheorem{remark}{Remark}

\theoremstyle{plain}
\newtoks\thehProclaim
\newtheorem*{Proclaim}{\the\thehProclaim}

\theoremstyle{definition}

\begin{document}

\title[Monotonicity formula for curve shortening flow]{Stabilization technique applied to curve shortening flow in the plane}

\author{Hayk Mikayelyan}

\address{Hayk Mikayelyan \\Mathematcal Sciences\\
Xi'an Jiaotong-Liverpool University\\
RenAiLu 111,
215123 Suzhou (SIP) \\
Jiangsu Prov., PR China}
              
\email{Hayk.Mikayelyan@xjtlu.edu.cn}


              
\begin{abstract}
The method proposed by T. I. Zelenjak is applied to the 
mean curvature flow in the plane. A new type of monotonicity 
formula for star-shaped curves is obtained.  
\end{abstract}

\date{1/DEC/2014}

\maketitle

\section{Introduction}

One of the classical problems combining geometry and PDEs is the mean curvature flow. 
Gerhard Huisken proved that convex surfaces converge in finite time to points in 
asymptotically spheric fashion (see  \cite{H1}). In dimension two this result was proven 
by M. Gage and R. Hamilton in \cite{GH}.  In \cite{G} Grayson showed that in the plane any closed 
embedded curve shrinks to a convex one in finite time and thus also shrinks to a point. This result is not true 
in higher dimensions where other types of singularities may occur if the initial curve is not convex (see \cite{Zh}).

A powerful tool in solving  is the monotonicity formula of Gerhard Huisken (see \cite{H2}), 
which plays a central role in many proofs.

In the present paper we apply a general method developed by T. Zelenjak in \cite{Z} to mean curvature 
flow in the plane, and derive the monotonicity formula of Huisken. We also derive a new monotonicity formula
for star-shaped curves. The presented approach is general and systematic, and we believe can be very useful 
in generalizations of the mean curvature flow, where no monotonicity formula is known. Our main 
motivation was the derivation of such a monotonicity formula for the anisotropic mean curvature flow, which still
remains a challenge.  

\section{The formulation of the problem}

We consider a closed curve in $\mathbb{R}^2$ moving by its curvature
with an anisotropy given by a function $g$:
$$
\partial_t \gamma = g(\nu)\kappa\nu,
$$
where $\gamma:\mathbb{R}_+ \times S^1\to \mathbb{R}^2$ is the curve parametrization,
$\kappa$ is the curvature and $\nu$ is the normal vector.

Note that in this form we fix a certain parametrization which has no
tangential component. For a general parametrization we will get
\begin{equation}\label{weaker}
\partial_t \gamma \cdot\nu= g(\nu)\kappa.
\end{equation}

If we take now $\gamma(t,x)=\left(\begin{array}{cc} u_1(t,x) \\u_2(t,x)
\end{array}\right)$ we get the following
$$
\left(
\begin{array}{cc} \partial_t u_1\\\partial_t u_2
\end{array}
\right)=
g(u_1',u_2')
\frac{-u_1''u_2'+u_1'u_2''}{({u_1'}^2+{u_2'}^2)^2}
\left(
\begin{array}{cc} -u_2' \\u_1'
\end{array}
\right)
$$
where $ ' $ means the $x$-derivative.

Assume the first singularity appears at point $0$ after finite time $T$.
We rescale the parametrization in the following way
$$
\tau=-\log(T-t), \tilde{\gamma}(\tau,x)=(T-t)^{-\frac{1}{2}}\gamma(t,x)
$$
and arrive at
\begin{equation}  \label{main}
\left(
\begin{array}{cc} \partial_\tau v_1\\ \partial_\tau v_2
\end{array}
\right)=
\frac{1}{2}
\left(
\begin{array}{cc} 
 v_1\\ v_2
 \end{array}
 \right)
 +g(v_1',v_2')
\frac{-v_1''v_2'+v_1'v_2''}{({v_1'}^2+{v_2'}^2)^2}
\left(
\begin{array}{cc} -v_2' \\v_1'
\end{array}
\right).
\end{equation}

In the paper we will carry out the main part of the computations for the 
anisotropic flow, but we are able to do the final part of the computations 
only for the isotropic case.
\begin{remark}
Note that in the classical case $g \equiv 1$ the stationary solution is the circle with the radius $\sqrt{2}$.
\end{remark}

\subsection*{Main result}
For the 
solutions of (\ref{main}),
which are star-shaped with respect to the origin, we prove the monotonicity 
formula
\begin{multline}
\frac{d}{d\tau}\int_{S^1} \sqrt{\frac{v'^2_1+v'^2_2}{v_1^2+v_2^2}} \Big(  f(\psi)+\Big( 
\frac{\log (v_1^2+v_2^2)}{2}+
 \frac{v_1^2+v_2^2}{4}
 \Big) \cos\psi  \Big) dx=\\
-\int_{S^1}|\partial_\tau \gamma\cdot\nu|^2  \sqrt{\frac{v'^2_1+v'^2_2}{v_1^2+v_2^2}}  \frac{1}{\cos \psi}dx.
\end{multline}
where $\psi$ is the angle 
between the outer normal direction $(v'_2,-v'_1)$ and the position vector $(v_1,v_2)$, and 
the function $f$
is defined in (\ref{ffff}), see also Figure \ref{fig-nc}.

\section{Monotonicity formula by Zelenjak's approach}

In this section we adapt the method proposed by T. I. Zelenjak in \cite{Z} to mean curvature flow. 

For this system we want to obtain a monotonicity formula of the form
\begin{equation}\label{monform}
\frac{d}{d\tau}\int_{S^1} F(v_1,v_2,v_1',v_2')dx=
-\int_{S^1}|\partial_\tau \gamma\cdot\nu|^2\rho(v_1,v_2,v_1',v_2')dx,
\end{equation}
where $\rho$ is positive. 

Note that in the classical case $g\equiv 1$ the well known Huisken's monotonicity 
formula (see \cite{H2},\cite{Zh}) in this
notations will correspond to the 
$$
F(\xi_1,\xi_2,\eta_1,\eta_2)=
\rho(\xi_1,\xi_2,\eta_1,\eta_2)=e^{-\frac{|\xi|^2}{4}}|\eta|.
$$
Diffentiating the left hand side of (\ref{monform}) and integrating by parts 
we get 
\begin{align}\label{tau1}
\partial_\tau v_1\left[\frac{\partial F}{\partial\xi_1}-
\frac{\partial^2 F}{\partial\xi_1\partial\eta_1}v_1'-
\frac{\partial^2 F}{\partial\xi_2\partial\eta_1}v_2'-
\frac{\partial^2 F}{\partial\eta_1^2}v_1''-
\frac{\partial^2 F}{\partial\eta_1\partial\eta_2}v_2''
\right]+\\
\label{tau2}
\partial_\tau v_2\left[\frac{\partial F}{\partial\xi_2}-
\frac{\partial^2 F}{\partial\xi_1\partial\eta_2}v_1'-
\frac{\partial^2 F}{\partial\xi_2\partial\eta_2}v_2'-
\frac{\partial^2 F}{\partial\eta_1\partial\eta_2}v_1''-
\frac{\partial^2 F}{\partial\eta_2^2}v_2''
\right].
\end{align}
In the right hand side of (\ref{monform}) using (\ref{main}) we obtain
\begin{align}
-\rho(\xi,\eta)
\frac{-\partial_\tau v_1 v_2'+\partial_\tau v_2 v_1'}{({v_1'}^2+{v_2'}^2)^\frac{1}{2}}
\left(\frac{-v_1v_2'+v_2v_1'}{2({v_1'}^2+{v_2'}^2)^\frac{1}{2}}+
g(v_1',v_2')\frac{-v_1''v_2'+v_2''v_1'}{({v_1'}^2+{v_2'}^2)^\frac{3}{2}}
\right)=\\
\label{rhotau1}
-\rho\partial_\tau v_1\left[
\frac{v_1{v_2'}^2-v_2v_1'v_2'}{2({v_1'}^2+{v_2'}^2)}+
g(v_1',v_2')\frac{{v_2'}^2}{({v_1'}^2+{v_2'}^2)^2}v_1''-
g(v_1',v_2')\frac{v_1'v_2'}{({v_1'}^2+{v_2'}^2)^2}v_2''
\right]\\
\label{rhotau2}
-\rho\partial_\tau v_2\left[
\frac{v_2{v_1'}^2-v_1v_1'v_2'}{2({v_1'}^2+{v_2'}^2)}-
g(v_1',v_2')\frac{v_1'v_2'}{({v_1'}^2+{v_2'}^2)^2}v_1''+
g(v_1',v_2')\frac{{v_1'}^2}{({v_1'}^2+{v_2'}^2)^2}v_2''
\right].
\end{align}
\begin{remark}
Note that we do not use (\ref{main}) but its weak form 
similar to (\ref{weaker}), which
means that what we get will work for any parametrization of the curve.
\end{remark}

We now require that the square brackets of (\ref{tau1}) and (\ref{rhotau1})
as well as (\ref{tau2}) and (\ref{rhotau2}) be equal. Moreover, we
require that
\begin{multline}\label{fandrho}
D^2_\eta F(\xi,\eta)=\rho(\xi,\eta)g(\eta)
\left(\begin{array}{cc} 
\frac{\eta_2^2}{(\eta_1^2+\eta_2^2)^2} & -\frac{\eta_1\eta_2}{(\eta_1^2+\eta_2^2)^2}\\ 
-\frac{\eta_1\eta_2}{(\eta_1^2+\eta_2^2)^2} & \frac{\eta_1^2}{(\eta_1^2+\eta_2^2)^2}
\end{array}\right)\\ =
\rho(\xi,\eta)g(\eta)|\eta|^{-1}D^2|\eta|
\end{multline}
where $g$ is homogeneous of order $0$. 

Introducing radial coordinates $(|\eta|, \phi)$ for $\eta$  
it is easy to check that for a given $\rho$
one can find an $F$ satisfying (\ref{fandrho}) if and only if
$$
\rho(\xi,\eta)=c(\xi,\eta)|\eta|,
$$ 
where $c$ is homogeneous of order $0$ with respect to $\eta$, and
$$
\int_0^{2\pi} c(\xi,\phi)g(\phi)\cos\phi d\phi=
\int_0^{2\pi} c(\xi,\phi)g(\phi)\sin\phi d\phi=0,\,\,\,\text{for all $\xi$}.
$$
Moreover, $F$ is homogeneous of order 1 in $\eta$ variable
and we can write $F(\xi, \eta)=f(\xi,\phi)|\eta|$. The formula
(\ref{fandrho}) becomes now
\begin{equation}\label{fandrho1}
D^2 F=(\partial_{\phi\phi}f+f)D^2|\eta|=c(\xi,\phi)g(\phi)D^2|\eta|.
\end{equation}
The solution of $f''+f=h$ can be calculated by the following formula
\begin{equation}\label{intconv}
f(\phi)=c_1 \cos \phi + c_2 \sin \phi + \int_0^\phi h(\tau)\sin(\phi-\tau)d\tau.
\end{equation}



What now remains to achieve our aim (\ref{tau1})=(\ref{rhotau1}) and (\ref{tau2})=(\ref{rhotau2}), is to make sure that 
\begin{equation}\label{remains1}
\frac{\partial F}{\partial\xi_1}-
\frac{\partial^2 F}{\partial\xi_1\partial\eta_1}\eta_1-
\frac{\partial^2 F}{\partial\xi_2\partial\eta_1}\eta_2=
\rho\frac{-\xi_1\eta_2^2+\xi_2\eta_1\eta_2}{2(\eta_1^2+\eta_2^2)}
\end{equation}
and 
\begin{equation}\label{remains2}
\frac{\partial F}{\partial\xi_2}-
\frac{\partial^2 F}{\partial\xi_1\partial\eta_2}\eta_1-
\frac{\partial^2 F}{\partial\xi_2\partial\eta_2}\eta_2=
\rho\frac{-\xi_2\eta_1^2+\xi_1\eta_1\eta_2}{2(\eta_1^2+\eta_2^2)}.
\end{equation}
After differentiating this equations in $\eta_1$ and $\eta_2$ respectively we obtain
$$
-\frac{\partial^3 F}{\partial\xi_1\partial\eta_1^2}\eta_1-
\frac{\partial^3 F}{\partial\xi_2\partial\eta_1^2}\eta_2=\partial_{\eta_1}
\left(
\rho\frac{-\xi_1\eta_2^2+\xi_2\eta_1\eta_2}{2(\eta_1^2+\eta_2^2)}\right)
$$
and
$$
-\frac{\partial^3 F}{\partial\xi_1\partial\eta_2^2}\eta_1-
\frac{\partial^3 F}{\partial\xi_2\partial\eta_2^2}\eta_2=
\partial_{\eta_2}
\left(
\rho\frac{-\xi_2\eta_1^2+\xi_1\eta_1\eta_2}{2(\eta_1^2+\eta_2^2)}\right),
$$
where we can substitute the value of $D_\eta^2 F$ from 
(\ref{fandrho}). The two equations turn out to be the same and can be written in terms of $c$ as follows
\begin{equation}\label{mainC}
2g(\eta)\langle \eta,D_\xi c\rangle-|\eta|^2\langle\xi,D_\eta c\rangle=-c\cdot\langle\xi,\eta\rangle.
\end{equation}
\begin{remark}\label{remdifeta12}
If we differentiate (\ref{remains1}) with respect to $\eta_2$ and add (\ref{remains2})
differentiated with respect to $\eta_1$ we will get the same (\ref{mainC}).
\end{remark}

Taking $c=e^b$ and rewriting (\ref{mainC}) in polar coordinates in $\eta$ variable we arrive at 
$$
\nabla_{\xi_1,\xi_2,\phi} b \cdot \left(\begin{array}{ccc} 2g(\phi)\cos \phi \\
2g(\phi) \sin \phi \\
\xi_1 \sin \phi - \xi_2 \cos \phi\end{array}\right)=-\xi_1\cos\phi
-\xi_2\sin\phi.
$$
The coordinate transformation 
$$
\begin{array}{lll} \tilde{\xi}_1=\xi_1 \cos \phi +\xi_2\sin \phi\\
\tilde{\xi}_2= \xi_1\sin \phi -\xi_2 \cos \phi\\
\tilde{\phi}=\phi
\end{array},
$$
brings us to the following first order linear PDE 
\begin{equation}\label{1-PDE}
\nabla_{ \tilde{\xi}_1 \tilde{\xi}_2,\phi} b \cdot \left(\begin{array}{ccc} 2g(\phi)- \tilde{\xi}_2^2
\\
 \tilde{\xi}_1 \tilde{\xi}_2\\
 \tilde{\xi}_2\end{array}\right)=- \tilde{\xi}_1.
\end{equation}

In the classical case $g\equiv 1$ the solution  
$b( \tilde{\xi},\phi)=-\frac{| \tilde{\xi}|^2}{4}$ gives us Huisken's famous monotonicity formula. 

It remains a 
challenge to find a solution to (\ref{1-PDE}) in the general anisotropic case which would correspond to the Huisken's one.


\section{A new monotonicity formula in isotropic case}


Let us observe that another obvious 
solution to (\ref{1-PDE}) different from Huisken's one is
\begin{equation}\label{logsol}
b( \tilde{\xi},\phi)=-\log | \tilde{\xi}_2|.
\end{equation}
This gives the function
$$
\rho(\xi,\eta)=|\eta|c(\xi,\eta)= \frac{|\eta|}{|\xi_1\sin\phi -\xi_2\cos\phi|}=\frac{|\eta|^2}{|\langle \xi , \eta_\nu  \rangle|},
$$
where $\eta_\nu$ is the vector $\eta$ rotated by 90 degrees clockwise and thus showing in outer normal 
direction.

Obviously we cannot solve (\ref{fandrho1}) globally because $\rho$ is not integrable, but if we assume 
that our domain is always star-shaped with respect to the origin and the angle $\psi$ 
between $\xi$ and $\eta_\nu$
will remain between $-\pi/2$ and $\pi/2$, we can solve (\ref{fandrho1}) locally. Now we 
just solve the equation  
\begin{equation}
\label{feqeq}
\partial_{\psi\psi}f+f=\frac{1}{\cos\psi}
\end{equation}
in the interval $(-\pi/2,\pi/2)$.
The general solution is 
\begin{equation}
\label{fgen}
f(\psi)+a(\xi)\cos\psi+b(\xi)\sin\psi,
\end{equation}
where 
\begin{equation}
\label{ffff}
f(\psi)=\psi\sin\psi+\cos\psi\log(\cos\psi).
\end{equation}
Let us first take $a=b=0$.

\begin{center}

\includegraphics[scale=0.6]{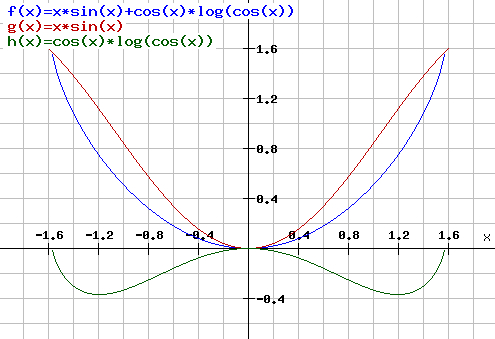}\label{fig-nc}

{\bf Figure \ref{fig-nc}}

\end{center}

As one can see from the graph that $f(\psi)$ is a positive, bounded, convex, even function in the 
interval $(-\pi/2,\pi/2)$.
The corresponding function $F$ is 
$$
F(\xi,\eta)=\frac{|\eta|}{|\xi|}f(\psi)=\frac{|\eta|}{|\xi|}(\psi\sin\psi+\cos\psi\log(\cos\psi)),
$$
where $\psi$ is the angle between position vector $\xi=(v_1,v_2)$ and the outer normal $\nu$
showing in the direction $(\eta_2,-\eta_1)=(v'_2,-v'_1)$.

Now we need to check whether the function $F$ satisfies the equations (\ref{remains1}) and 
(\ref{remains2}). The answer is no. We obtain in (\ref{remains1}) and 
(\ref{remains2})
$$
\frac{\eta_2}{|\xi|^2}\not= -\frac{\eta_2}{2}
$$
and 
$$
-\frac{\eta_1}{|\xi|^2}\not= \frac{\eta_1}{2}
$$
respectively (see Remark \ref{remdifeta12}).
This means we have an additional term in the formula (\ref{monform})
\begin{multline}\label{monformadd}
\frac{d}{d\tau}\int_{S^1} \sqrt{\frac{v'^2_1+v'^2_2}{v_1^2+v_2^2}}   f(\psi) dx+
\int_{S^1}|\partial_\tau \gamma\cdot\nu|^2  \sqrt{\frac{v'^2_1+v'^2_2}{v_1^2+v_2^2}}  \frac{1}{\cos \psi}dx=\\
-\int_{S^1} ( v'_2\partial_\tau v_1- v'_1\partial_\tau v_2)\Big(\frac{1}{2}+\frac{1}{v_1^2+v_2^2}  \Big)dx.
\end{multline}


\section{The ``repaired'' formula}\label{secrep}
In order to obtain a monotonicity formula without additional terms we need to go back to the 
general solution of  (\ref{feqeq}). The idea is that  
by adding a term linear in $\eta$ to $F$ we do not create problems in (\ref{fandrho}), so 
let us find a function $a(r)$
such that the function
$$
F(\xi,\eta)=\frac{|\eta|}{|\xi|}f(\psi)+  a(|\xi|) |\eta| \cos \psi
$$
solves (\ref{remains1}) and 
(\ref{remains2}), so we do not have additional terms. 
Substituting $F$ we obtain
$$
\frac{\eta_2}{|\xi|^2}  -  \eta_2\Big( a'(|\xi|) +\frac{a(|\xi|)}{|\xi|}  \Big)
=- \frac{\eta_2}{2}
$$
and 
$$
-\frac{\eta_1}{|\xi|^2}+ \eta_1\Big( a'(|\xi|) +\frac{a(|\xi|)}{|\xi|}  \Big)  
= \frac{\eta_1}{2}
$$
respectively, and now need to solve
\begin{equation}
ra'(r)+a(r)=\frac{r}{2}+\frac{1}{r}.
\end{equation}
The solution is $a(r)=\frac{r}{4}+\frac{\log r}{r}$ and 
$$
F(\xi,\eta)=\frac{|\eta|}{|\xi|}f(\psi)+ |\eta| \Big(\frac{|\xi|}{4}+\frac{\log |\xi|}{|\xi|} \Big)\cos \psi.
$$
Thus we obtain the following monotonicity formula
\begin{multline}\label{monforrr}
\frac{d}{d\tau}\int_{S^1} \sqrt{\frac{v'^2_1+v'^2_2}{v_1^2+v_2^2}} \Big(  f(\psi)+\Big( 
\frac{\log (v_1^2+v_2^2)}{2}+
 \frac{v_1^2+v_2^2}{4}
 \Big) \cos\psi  \Big) dx=\\
-\int_{S^1}|\partial_\tau \gamma\cdot\nu|^2  \sqrt{\frac{v'^2_1+v'^2_2}{v_1^2+v_2^2}}  \frac{1}{\cos \psi}dx.
\end{multline}



\subsection*{Acknowledgment} The author is is grateful to Georg Weiss and Hyunsuk Kang for inspiring discussions.








\end{document}